\documentclass[runningheads]{llncs}
\usepackage{bdlnat}
\usepackage[hidelinks]{hyperref}
\hypersetup{colorlinks=true,linkcolor=blue,citecolor=blue,urlcolor=blue}
\usepackage{orcidlink}
\usepackage{version}

\newcommand{\doi}[1]{{doi:\href{https://doi.org/#1}{\nolinkurl{#1}}}}
\renewcommand{\url}[1]{{\href{#1}{\nolinkurl{#1}}}}

\spnewtheorem{sdef}{Definition}{\bfseries}{\rmfamily}

\title{The Most Natural Paradefinite Logic \\ Relative to Classical Logic}
\author{C.A. Middelburg\,\orcidlink{0000-0002-8725-0197}}
\institute
 {Informatics Institute, Faculty of Science, University of Amsterdam \\
  Science Park~900, 1098~XH Amsterdam, the Netherlands \\
  \href{mailto:C.A.Middelburg@uva.nl}{C.A.Middelburg@uva.nl}}

\titlerunning
 {The Most Natural Paradefinite Logic Relative to Classical Logic}
\authorrunning
 {C.A. Middelburg}

\begin{document}

\maketitle

\begin{abstract}
A paradefinite logic is a logic that can serve as the underlying logic 
for theories that are inconsistent or incomplete.
A well-known paradefinite logic is Belnap-Dunn logic.
Various expansions of Belnap-Dunn logic have been studied in the 
literature.
In this note, it is argued that the most natural paradefinite logic 
relative to classical logic is the expansion of Belnap-Dunn logic with 
a falsity connective and an implication connective for which the 
standard deduction theorem holds.
\begin{msc-classcode} 
03B50 (Primary)\,\, 03B05, 03B53 (Secondary) 
\end{msc-classcode}
\end{abstract}

\section{Introduction}
\label{INTRO}

Belnap-Dunn logic (BD)~\cite{Bel77a,Bel77b}, also called the logic 
of First-Degree Entailment, is known for the fact that it can serve as 
the underlying logic of theories that are inconsistent or incomplete.
BD is closely related to three other well-known logics: its logical 
consequence relation is included in the logical consequence relations of 
Priest's Logic of Paradox (\LP)~\cite{Pri79a}, Kleene's strong 
$3$-valued logic (\Kiii)~\cite{Kle52a}, and the version of classical 
propositional logic with the same connectives as BD.
BD can serve as the underlying logic of theories that are inconsistent 
or incomplete because it is paraconsistent and paracomplete.
That is, like \LP, BD is a logic in which not every formula is a logical 
consequence of each two formulas of which one is the negation of the 
other and, like \Kiii, BD is a logic in which not, for each two formulas 
of which one is the negation of the other, one or the other is a logical 
consequence of every set of formulas.

Expansions of BD (or a first-order extension of it) that are studied in 
earlier papers are usually expansions with one or more connectives that 
are not known from classical propositional logic.
Examples are \BDD~\cite{SO14a}, \EFDEN~\cite{CC20a}, 
\FivCC~\cite{KZ20a}, and \QLETF~\cite{ARCC22a}.  
This note concerns an expansion of BD with connectives known from 
classical propositional logic, namely a falsity connective and an 
implication connective for which the standard deduction theorem holds.
This expansion will be referred to as $\BDL$ in this note.
It has been treated in several earlier papers, 
including~\cite{AA96a,AA98a,AA17a,Pyn99a}, but without exception quite 
casually.

In most of this note, the phrase `a connective known from classical 
logic' is used rather loosely.
Roughly speaking, a connective of a non-classical logic is considered a 
connective known from classical logic if it has the same properties with 
respect to logical consequence as a connective of some version of 
classical logic.

Following~\cite{AA17a}, among others, logics that can serve as the 
underlying logic of theories that are inconsistent or incomplete are 
referred to in this note as paradefinite logics.

The structure of this note is as follows.
First, preliminaries concerning propositional logics are presented (Section~\ref{PRELIMINARIES}).
Next, it is argued that the most natural paradefinite logic relative to 
classical logic is the expansion of Belnap-Dunn logic with a falsity 
connective and an implication connective for which the standard 
deduction theorem holds (Section~\ref{NATURAL}). 
Finally, some concluding remarks are made (Section~\ref{CONCLUSIONS}).

\section{Propositional Logics}
\label{PRELIMINARIES}

The scope of this note is limited to propositional logics equipped with 
a structural and non-trivial Scott consequence relation.
Precise definitions are required for various notions relevant to logics 
of the kind considered in this note.
The relevant definitions are presented in this section.

The language of a propositional logic is defined by way of an alphabet 
that consists of propositional variables and logical connectives.
\begin{sdef}
An \emph{alphabet} of a language of a propositional logic is a couple
$\LAlph = \tup{\PVar,\indfam{\Conn{n}}{n \in \Nat}}$, where:
\begin{itemize}
\item
$\PVar$ is a countably infinite set of \emph{propositional variables};
\item
$\indfam{\Conn{n}}{n \in \Nat}$ is an $\Nat$-indexed family of pairwise 
disjoint sets;
\item
for each $n \in \Nat$, $\Conn{n}$ is a finite set of 
\emph{logical connectives of arity $n$};
\item
$\Union \set{\Conn{n} \where n \in \Nat}$ is a non-empty finite set.
\end{itemize}
\end{sdef}

The language over an alphabet consists of formulas.
They are constructed according to the formation rules given below.

\begin{sdef}
Let $\LAlph = \tup{\PVar,\indfam{\Conn{n}}{n \in \Nat}}$ be an alphabet.
Then the set $\LForm^\LAlph$ of all \emph{formulas} over $\LAlph$ is 
inductively defined by the following formation rules:
\begin{itemize}
\item
if $p \in \PVar$, then $p \in \LForm^\LAlph$;
\item
if ${\Diamond} \in \Conn{0}$, then ${\Diamond} \in \LForm^\LAlph$;
\item
if ${\Diamond} \in \Conn{n+1}$ and
$A_1,\ldots,A_{n+1} \in \LForm^\LAlph$, then 
${\Diamond}(A_1,\ldots,A_{n+1}) \in \LForm^\LAlph$.
\end{itemize}
\pagebreak[2]
The set of all \emph{atomic formulas} over $\LAlph$, written 
$\LAForm^\LAlph$, is the set $\PVar$ of propositional variables.
\end{sdef}
\pagebreak[2]
The following will sometimes be used without mentioning (with or without 
decoration):
$p$ and $q$ as meta-variables ranging over all propositional variables 
from $\PVar^\LAlph$,
$A$ and $B$ as meta-variables ranging over all formulas from 
$\LForm^\LAlph$, and
$\Gamma$ and $\Delta$ as meta-variables ranging over all sets of 
formulas from $\LForm^\LAlph$.

\begin{sdef}
Let $\LAlph = \tup{\PVar,\indfam{\Conn{n}}{n \in \Nat}}$ be an alphabet.
Then an \emph{$\LAlph$-substi\-tution} of formulas from $\LForm^\LAlph$ 
for variables from $\PVar$ is a function 
$\funct{\sigma}{\PVar}{\LForm^\LAlph}$.
An $\LAlph$-substitution $\sigma$ extends to the function 
$\funct{\sigma^*}{\LForm^\LAlph}{\LForm^\LAlph}$ that is 
recursively defined as follows:
\begin{ldispl}
\begin{array}[t]{r@{\;}c@{\;}l}
\sigma^*(p)                                & = & \sigma(p),\; \\
\sigma^*({\Box})                           & = & {\Box},\; \\
\sigma^*({\Diamond}(A_1, \ldots, A_{n+1})) & = &
 {\Diamond}(\sigma^*(A_1), \ldots,\sigma^*(A_{n+1})),
\end{array}
\end{ldispl}%
for ${\Box} \in \Conn{0}$ and ${\Diamond} \in \Conn{n+1}$.
\end{sdef}
We write $\sigma(A)$ for $\sigma^*(A)$ and
$\sigma(\Gamma)$ for $\set{\sigma^*(A) \where A \in \Gamma}$. 

We use the notational conventions to write $(\Diamond A)$ instead of 
$\Diamond(A)$ and $(A_1 \Diamond A_2)$ instead of $\Diamond(A_1,A_2)$ 
and to omit parenthesis where it does not lead to syntactic ambiguities 
if the previous convention is used.

\begin{sdef}
Let $\LAlph = \tup{\PVar,\indfam{\Conn{n}}{n \in \Nat}}$ be an alphabet.
Then an 
\emph{(ordinary) logical consequence relation for $\LForm^\LAlph$} is a 
binary relation $\LCon$ on $\pset(\LForm^\LAlph)$ that satisfies the 
following conditions:%
\footnote
{As usual, we write $\Gamma,\Gamma'$ for $\Gamma \union \Gamma'$ and $A$
 for $\set{A}$.}
\begin{itemize}
\item
if $\Gamma \inter \Delta \neq \emptyset$ then $\Gamma \LCon \Delta$;
\item
if $\Gamma \LCon \Delta$, $\Gamma \subseteq \Gamma'$, and 
$\Delta \subseteq \Delta'$ then $\Gamma' \LCon \Delta'$;
\item
if $\Gamma \LCon \Delta,A$ and $A,\Gamma' \LCon \Delta'$ then 
$\Gamma,\Gamma' \LCon \Delta,\Delta'$.
\end{itemize}
\end{sdef}
The qualification \emph{ordinary} in the above definition emphasizes 
that other kinds of logical consequence relation are considered in the 
literature on logic.
In this note, only ordinary logical consequence relations are 
considered.
Therefore, in the sequel, an ordinary logical consequence relation is 
simply called a logical consequence relation.
Ordinary logical consequence relations are also known as Scott 
consequence relations.

In this note, only propositional logics whose logical consequence 
relation satisfies the condition of being structural and the condition 
of being non-trivial are considered.
\begin{sdef}
Let $\LAlph = \tup{\PVar,\indfam{\Conn{n}}{n \in \Nat}}$ be an alphabet,
and let $\LCon$ be a logical consequence relation for $\LForm^\LAlph$.
Then $\LCon$ is \emph{structural} if it satisfies the following 
condition:
\begin{itemize}
\item[]
if $\Gamma \LCon \Delta$ and $\sigma$ is an $\LAlph$-substitution then 
$\sigma(\Gamma) \LCon \sigma(\Delta)$.
\end{itemize}
\end{sdef}
\begin{sdef}
Let $\LAlph = \tup{\PVar,\indfam{\Conn{n}}{n \in \Nat}}$ be an alphabet,
and let $\LCon$ be a logical consequence relation for $\LForm^\LAlph$.
Then $\LCon$ is \emph{non-trivial} if it satisfies the following 
condition:
\begin{itemize}
\item[]
there exist non-empty $\Gamma$ and $\Delta$ such that not 
$\Gamma \LCon \Delta$.
\end{itemize}
\end{sdef}

\begin{sdef}
A \emph{(propositional) logic} is a couple $(\LAlph,\LCon)$, where:
\begin{itemize}
\item
$\LAlph$ is an alphabet;
\item
$\LCon$ is a structural and non-trivial logical consequence relation for 
$\LForm^\LAlph$.
\end{itemize}
\end{sdef}
In this note, only propositional logics are considered.
Therefore, in the sequel, a propositional logic is simply called a 
logic.

The condition of being structural is usually required in the 
definitions of propositional logics.
The condition of being non-trivial is often not required in the 
definitions of propositional logics, but is convenient for excluding
trivial logics.

\section{Naturally from Classical Logic to a Paradefinite Logic}
\label{NATURAL}

The particular version of classical logic  and the particular expansion 
of Belnap-Dunn logic considered in this note are called \CL\ and \BDL, 
respectively.

The language of \CL\ and the language of \BDL\ are the same.
This means that \CL\ and \BDL\ have a common alphabet.
\begin{sdef}
\label{def-LAlph-CL}
The alphabet $\CLAlph$ of the language of an instance of \CL\ or \BDL\ 
is a couple $(\PVar,\indfam{\Conn{n}}{n \in \Nat})$, where:
\begin{itemize}
\item
$\PVar$ is a countably infinite set of propositional variables;
\item
$\Conn{0} = \set{\False}$;
\item
$\Conn{1} = \set{\Not}$;
\item
$\Conn{2} = \set{\CAnd,\COr,\IImpl}$; 
\item
$\Conn{n+3} = \emptyset$ for each $n \in \Nat$.
\end{itemize}
Each choice of $\PVar$ gives rise to a different instance of 
\CL\ and \BDL. 
In this note, a fixed but arbitrary choice of $\PVar$ is assumed.
\end{sdef}

Arguments for the choice of connectives are:
\begin{itemize}
\item
any expansion of BD must include the connectives $\CAnd$, $\COr$, and 
$\Not$ because these are the connectives of BD;
\item
the expansion of BD with both the connectives $\IImpl$ and $\False$ has 
greater expressive power than an expansion of BD with only one of them;
\item
an expansion of BD with the connectives $\IImpl$ and $\False$ and other 
connectives known from classical logic does not have more expressive 
power than the expansion of BD with only the connectives $\IImpl$ and 
$\False$;
\item
an expansion of BD with connectives not known from classical logic 
does not deserve to be qualified as the most natural paradefinite logic 
relative to classical logic. 
\end{itemize}
Moreover, this choice of connectives yields a suitable language for 
the most natural paradefinite logic relative to classical logic: 
it guarantees that, for any connective available or definable in the 
version of classical logic, a connective with the same properties with 
respect to logical consequence is available or definable in the 
paradefinite logic.

It is worth mentioning here that, although the connectives of BD are 
$\CAnd$, $\COr$, and~$\Not$, the implication connective $\IImpl$ for 
which the standard inference theorem holds and the falsity connective 
$\False$ are not definable in BD.
A relatively unknown consequence of expanding BD with the connectives 
$\IImpl$ and $\False$ is that several interesting connectives not known 
from classical logic become definable (see~\cite{Mid24b}, Section~6).
Additional connectives not known from classical logic are needed to 
obtain an expansion of BD with more expressive power than the 
expansion of BD with the connectives $\IImpl$ and $\False$.

Henceforth, we write $\LForm$ and $\LAForm$ instead of $\LForm^\CLAlph$ 
and $\LAForm^\CLAlph$, respectively.

The logical consequence relation of a logic is usually explicitly 
defined using a logical matrix.
In this note, we give implicit definitions of the logical consequence 
relations of \CL\ and \BDL\ instead. 
In either case, the implicit definition shows more directly how the 
logical consequence relation concerned and the different connectives are 
related.
\begin{sdef}
\label{def-LCon-CL}
The logical consequence relation $\LCon_\CL$ of \CL\ is the smallest 
logical consequence relation such that 
for all $\Gamma, \Delta \subseteq \LForm$, and $A_1, A_2 \in \LForm$:
\[
\renewcommand{\arraystretch}{1.25}
\begin{array}[t]{r@{\;}c@{\;}l@{\;\;}l}
\False, \Gamma \LCon_\CL \Delta & \!\!\!\!,
\\
\Gamma \LCon_\CL \Delta, \Not A_1         & \mathrm{iff} &
\multicolumn{2}{l}
{A_1, \Gamma \LCon_\CL \Delta\;,}
\\
\Gamma \LCon_\CL \Delta, A_1 \CAnd A_2    & \mathrm{iff} &
\multicolumn{2}{l}
{\Gamma \LCon_\CL \Delta, A_1 \;\mathrm{and}\;
 \Gamma \LCon_\CL \Delta, A_2\;,}
\\
A_1 \COr A_2, \Gamma \LCon_\CL \Delta     & \mathrm{iff} &
\multicolumn{2}{l}
{A_1, \Gamma \LCon_\CL \Delta \;\mathrm{and}\;
 A_2, \Gamma \LCon_\CL \Delta\;,}
\\
\Gamma \LCon_\CL \Delta, A_1 \IImpl A_2   & \mathrm{iff} & 
A_1, \Gamma \LCon_\CL \Delta, A_2\;.
\end{array}
\]
\end{sdef}
The standard sequent calculus proof system for \CL\ is sound and 
complete with respect to the logical consequence relation defined using 
the logical matrix of \CL.
Because the rules of that proof system are invertible, it follows 
immediately that $\LCon_\CL$ as defined above coincides with the logical 
consequence relation defined using the logical matrix of \CL\ 
(cf.~\cite{Avr91c,Avr91a}).

The condition
$\Gamma \LCon_\CL \Delta, \Not A_1 \;\mathrm{iff}\;
 A_1, \Gamma \LCon_\CL \Delta$ 
in the above definition of $\LCon_\CL$ is a general condition concerning 
negation.
It can be replaced by a number of more specific conditions concerning
negation.
\begin{theorem}
\label{theorem-LCon-CL}
The logical consequence relation $\LCon_\CL$ of \CL\ is the smallest 
logical consequence relation such that 
for all $\Gamma, \Delta \subseteq \LForm$, and $A_1, A_2 \in \LForm$:
\pagebreak[2]
\[
\renewcommand{\arraystretch}{1.25}
\begin{array}[t]{r@{\;}c@{\;}l@{\;\;}l}
\False, \Gamma \LCon_\CL \Delta & \!\!\!\!,
\\
\Gamma \LCon_\CL \Delta, A_1 \CAnd A_2    & \mathrm{iff} &
\multicolumn{2}{l}
{\Gamma \LCon_\CL \Delta, A_1 \;\mathrm{and}\;
 \Gamma \LCon_\CL \Delta, A_2\;,}
\\
A_1 \COr A_2, \Gamma \LCon_\CL \Delta     & \mathrm{iff} &
\multicolumn{2}{l}
{A_1, \Gamma \LCon_\CL \Delta \;\mathrm{and}\;
 A_2, \Gamma \LCon_\CL \Delta\;,}
\\
\Gamma \LCon_\CL \Delta, A_1 \IImpl A_2   & \mathrm{iff} & 
A_1, \Gamma \LCon_\CL \Delta, A_2\;,
\\
\Not A_1, A_1 \LCon_\CL \False                   & \mathrm{and} &
\multicolumn{2}{l}
{\Not \False \LCon_\CL A_1, \Not A_1\;,}
\\
\Gamma \LCon_\CL \Delta, \Not \False & \!\!\!\!, 
\\
\Not (\Not A_1), \Gamma \LCon_\CL \Delta         & \mathrm{iff} &
\multicolumn{2}{l}
{A_1, \Gamma \LCon_\CL \Delta\;,}
\\
\Not (A_1 \CAnd A_2), \Gamma \LCon_\CL \Delta    & \mathrm{iff} &
\multicolumn{2}{l}
{\Not A_1, \Gamma \LCon_\CL \Delta \;\mathrm{and}\;
 \Not A_2, \Gamma \LCon_\CL \Delta\;,}
\\
\Gamma \LCon_\CL \Delta, \Not (A_1 \COr A_2)     & \mathrm{iff} &
\multicolumn{2}{l}
{\Gamma \LCon_\CL \Delta, \Not A_1 \;\mathrm{and}\;
 \Gamma \LCon_\CL \Delta, \Not A_2\;,}
\\
\Not (A_1 \IImpl A_2), \Gamma \LCon_\CL \Delta   & \mathrm{iff} & 
{A_1, \Not A_2, \Gamma \LCon_\CL \Delta\;.}
\end{array}
\]
\end{theorem}
\begin{proof}
Consider the sequent calculus proof system for \CL\ obtained by replacing 
the pair of invertible rules from the standard sequent calculus proof 
system for \CL\ that corresponds to the general condition concerning 
negation by the pairs of invertible rules that correspond to the more 
specific conditions concerning negation (cf.~\cite{Avr91a}).
To prove the theorem, it is sufficient to show that the resulting proof 
system is sound and complete.
This is easily done by induction on the length of the proof of a 
sequent, using the fact that the standard sequent calculus proof system 
for \CL\ is sound and complete.
\qed
\end{proof}

Theorem~\ref{theorem-LCon-CL} provides an alternative definition of 
$\LCon_\CL$.
If we get rid of the conditions $\Not A_1, A_1 \LCon_\CL \False$ and 
$\Not \False \LCon_\CL A_1, \Not A_1$ from this alternative definition, 
we obtain the implicit definition of the logical consequence relation of 
\BDL.
\begin{sdef}
\label{def-LCon-BD}
The logical consequence relation $\LCon_\BDL$ of \BDL\ is the smallest 
logical consequence relation such that 
for all $\Gamma, \Delta \subseteq \LForm$, and $A_1, A_2 \in \LForm$:
\[
\renewcommand{\arraystretch}{1.25}
\begin{array}[t]{r@{\;}c@{\;}l@{\;\;}l}
\False, \Gamma \LCon_\BDL \Delta & \!\!\!\!,
\\
\Gamma \LCon_\BDL \Delta, A_1 \CAnd A_2    & \mathrm{iff} &
\multicolumn{2}{l}
{\Gamma \LCon_\BDL \Delta, A_1 \;\mathrm{and}\;
 \Gamma \LCon_\BDL \Delta, A_2\;,}
\\
A_1 \COr A_2, \Gamma \LCon_\BDL \Delta     & \mathrm{iff} &
\multicolumn{2}{l}
{A_1, \Gamma \LCon_\BDL \Delta \;\mathrm{and}\;
 A_2, \Gamma \LCon_\BDL \Delta\;,}
\\
\Gamma \LCon_\BDL \Delta, A_1 \IImpl A_2   & \mathrm{iff} & 
A_1, \Gamma \LCon_\BDL \Delta, A_2\;,
\\
\Gamma \LCon_\BDL \Delta, \Not \False & \!\!\!\!, 
\\
\Not (\Not A_1), \Gamma \LCon_\BDL \Delta         & \mathrm{iff} &
\multicolumn{2}{l}
{A_1, \Gamma \LCon_\BDL \Delta\;,}
\\
\Not (A_1 \CAnd A_2), \Gamma \LCon_\BDL \Delta    & \mathrm{iff} &
\multicolumn{2}{l}
{\Not A_1, \Gamma \LCon_\BDL \Delta \;\mathrm{and}\;
 \Not A_2, \Gamma \LCon_\BDL \Delta\;,}
\\
\Gamma \LCon_\BDL \Delta, \Not (A_1 \COr A_2)     & \mathrm{iff} &
\multicolumn{2}{l}
{\Gamma \LCon_\BDL \Delta, \Not A_1 \;\mathrm{and}\;
 \Gamma \LCon_\BDL \Delta, \Not A_2\;,}
\\
\Not (A_1 \IImpl A_2), \Gamma \LCon_\BDL \Delta   & \mathrm{iff} & 
{A_1, \Not A_2, \Gamma \LCon_\BDL \Delta\;.}
\end{array}
\]
\end{sdef}
The sequent calculus proof system for \BDL\ presented in~\cite{Mid24b} 
is sound and complete with respect to the logical consequence relation 
defined using the logical matrix of \BDL\ presented in~\cite{Mid24b}.
Because the rules of that proof system are invertible, it follows 
immediately that $\LCon_\BDL$ as defined above coincides with the logical 
consequence relation defined using that logical matrix
(cf.~\cite{Avr91c,Avr91a}).

The conditions $\Not A_1, A_1 \LCon_\CL \False$ and 
$\Not \False \LCon_\CL A_1, \Not A_1$ from the alternative definition
of $\LCon_\CL$ provided by Theorem~\ref{theorem-LCon-CL} represent the 
\emph{law of non-contradiction} (LNC) and the 
\emph{law of excluded middle} (LEM), respectively.
LNC is the only reason why \CL\ cannot serve as the underlying logic for
theories that are inconsistent and LEM is the only reason why \CL\ 
cannot serve as the underlying logic for theories that are incomplete.
Getting rid of LNC and LEM is all that is needed to obtain a logic that 
can serve as the underlying logic for theories that are inconsistent or 
incomplete.
\BDL\ can be thought as obtained in exactly this way.

Clearly, it is the case that ${\LCon_\CL} \diff {\LCon_\BDL}$ is
precisely the set of all classical logical consequences that exist due 
to either LNC or LEM.
The most inartificial paradefinite logic relative to \CL\ is the logic 
that differs from \CL\ only in that it lacks exactly those logical 
consequences.
This makes that \BDL\ deserves to be qualified as the most natural 
paradefinite logic relative to \CL.

\section{Concluding Remarks}
\label{CONCLUSIONS}

It has been argued in this note that \BDL\ is the most natural 
paradefinite logic relative to \CL.
It can further be argued along the same lines that (a) by getting rid of 
only the condition $\Not A_1, A_1 \LCon_\CL \False$ from the alternative 
definition of \CL, we obtain the implicit definition of the logical 
consequence relation of the most natural paraconsistent logic relative 
to classical logic (to wit \LP\ expanded with $\IImpl$ and $\False$) and 
(b) by getting rid of only the condition 
$\Not \False \LCon_\CL A_1, \Not A_1$ from the alternative definition of 
\CL, we obtain the implicit definition of the logical consequence 
relation of the most natural paracomplete logic relative to classical 
logic (to wit \Kiii\ expanded with $\IImpl$ and $\False$).

\sloppy
It is not difficult to see that getting rid of the conditions 
$\Not A_1, A_1 \LCon_\CL \False$ and 
$\Not \False \LCon_\CL A_1, \Not A_1$ from the alternative definition 
of $\LCon_\CL$ agrees with treating negations of propositional variables 
as additional propositional variables.
This allows for a very simple embedding of \BDL\ into \CL\ 
(see~\cite{Mid24b}, Theorem~12).

\bibliographystyle{splncs04}
\bibliography{PCL}

\end{document}